\documentclass[11pt]{article}

\newcommand{\eit}{{e^{i\theta}}}
\newcommand{\me}{{-1}}

\newcommand{\n}{\|}
\newcommand{\iy}{\infty}
\newcommand{\hal}{{(\al)}}
\newcommand{\hn}{{(n)}}

\newcommand{\al}{\alpha}

\newcommand{\eps}{\varepsilon}
\newcommand{\ga}{\gamma}

\newcommand{\de}{\delta}
\newcommand{\la}{\lambda}
\newcommand{\ph}{\varphi}

\newcommand{\tht}{\theta}

\newcommand{\om}{\omega}

\newcommand{\bC}{{\bf C}}

\newcommand{\bR}{{\bf R}}
\newcommand{\bT}{{\bf T}}

\newcommand{\cP}{{\mathcal P}}

\newcommand{\vsk}{\vspace{0mm}}
\newcommand{\vsg}{\vspace{3mm}}

\begin{document}

{\Large\bf
\begin{center}
{From} Toeplitz Eigenvalues through Green's Kernels to
Higher-Order Wirtinger-Sobolev Inequalities
\end{center}}

{ \bf
\begin{center}
A. B\"ottcher and H. Widom
\end{center}}

\begin{quote}
\renewcommand{\baselinestretch}{1.0}
\footnotesize
The paper is concerned with a sequence of constants which appear in
several problems. These problems include the minimal eigenvalue of certain
positive definite Toeplitz matrices, the minimal eigenvalue of some higher-order
ordinary differential operators, the norm of the Green kernels of these operators,
the best constant in a Wirtinger-Sobolev inequality, and the conditioning of a
special least squares problem. The main result of the paper gives the asymptotics
of this sequence.
\end{quote}

\section{{\large Introduction and main result}}\label{S1}

There is a sequence $c_1, c_2, c_3, \ldots$ of positive real numbers that emerges
in various contexts. Here are five of them.

\vsk
\noindent
{\bf Minimal eigenvalues of Toeplitz matrices.} Given a continuous function $a$ on the
complex unit circle $\bT$, we denote by $\{a_k\}_{k=-\iy}^\iy$ the sequence of the
Fourier coefficients,
\[a_k=\frac{1}{2\pi} \int_0^{2\pi} a(\eit)e^{-ik\tht}d\tht,\]
and by $T_n(a)$ the $n \times n$ Toeplitz matrix $(a_{j-k})_{j,k=1}^n$.
Suppose $a$ is of the form $a(t)=|1-t|^{2\al}b(t)$ ($t \in \bT$) where
$\al$ is a natural number and $b$ is a positive function on $\bT$ whose
Fourier coefficients are subject to the condition
$\sum_{k=-\iy}^\iy |k|\,|b_k| < \iy$. Then the matrix
$T_n(a)$ is positive definite and its smallest eigenvalue
$\la_{\min}(T_n(a))$ satisfies
\begin{equation}
\la_{\min}(T_n(a)) \sim \frac{c_\al}{n^{2\al}}\,b(1) \quad \mbox{as} \quad n \to \iy
\label{1.1}
\end{equation}
with a certain constant $c_\al \in (0,\iy)$ independent of $b$. Here and in what
follows $x_n \sim y_n$ means that $x_n/y_n \to 1$. Kac, Murdock, and Szeg\"o
\cite{KMS} proved that $c_1=\pi^2$, and Parter \cite{Par99} showed that
$c_2=500.5467$.

\vsk
\noindent
{\bf Minimal eigenvalues of differential operators.} For a natural number $\al$,
consider the boundary value problem
\begin{eqnarray}
& & \hspace{-10mm} (-1)^\al \,u^{(2\al)}(x)=v(x) \quad \mbox{for}
\quad x \in [0,1], \label{1.2}\\
& & \hspace{-10mm} u(0)=u'(0)=\ldots=u^{(\al-1)}(0)=0, \quad
u(1)=u'(1)=\ldots=u^{(\al-1)}(1)=0. \label{1.3}
\end{eqnarray}
The minimal eigenvalue of this boundary value problem can be shown to be just
$c_\al$. If $\al=3$, then the equation $-u^{(6)}=\la u$ is satisfied by
\[u(x)=\sum_{k=0}^5 A_k \exp \left(x\,\sqrt[6]{\la}\,
\exp\left(\frac{(2k+1)\pi i}{6}\right)\right),\]
and the $A_k$'s are the solution of a homogeneous linear $6 \times 6$
system with a matrix depending on $\la$. We found numerically that
the smallest $\la >0$ for which the determinant of this
matrix is zero is approximately $\la=61529$. Thus, $c_3=61529$.

\vsk
\noindent
{\bf Norms of Green's kernels.} Let $G_\al(x,y)$ be the Green kernel of problem (\ref{1.2}),
(\ref{1.3}). The solution to (\ref{1.2}), (\ref{1.3}) is then given by
\begin{equation}
u(x)=\int_0^1G_\al(x,y)v(y)dy. \label{1.4}
\end{equation}
It can be shown that $G_\al(x,y)$ is symmetric about the point
$(\frac{1}{2},\frac{1}{2})$ and that
\begin{equation}
G_\al(x,y)=\frac{x^\al y^\al}{[(\al-1)!]^2}\int_{\max(x,y)}^1
\frac{(t-x)^{\al-1}(t-y)^{\al-1}}{t^{2\al}}\,dt \label{1.5}
\end{equation}
for $x+y \ge 1$. Let $K_\al$ denote the integral operator defined by (\ref{1.4}).
It is clear that the minimal eigenvalue of (\ref{1.2}), (\ref{1.3}) equals the inverse
of the maximal eigenvalue of the (compact and positive definite) operator $K_\al$
on $L^2(0,1)$. As the maximal eigenvalue of $K_\al$ is its norm, we arrive at the
equality $1/c_\al=\n K_\al \n$.

\vsk
\noindent
{\bf Best constants in Wirtinger-Sobolev inequalities.} By a Wirtinger-Sobolev
inequality one means an inequality of the form
\begin{equation}
\int_0^1|u(x)|^2dx \le C \int_0^1|u^\hal (x)|^2dx, \label{1.6}
\end{equation}
where $u$ is required to satisfy certain additional (for example, boundary)
conditions. It is well known that the best constant $C$ for which (\ref{1.6})
is true for all $u \in C^\al[0,1]$ satisfying $\int_0^1 u(x)dx=0$ and $u^{(j)}(0)
=u^{(j)}(1)$ for $0 \le j \le \al-1$ is equal to $1/(2\pi)^{2\al}$. However,
problem (\ref{1.2}), (\ref{1.3}) leads to (\ref{1.6}) with the additional constraints
(\ref{1.3}). In this case the best constant in (\ref{1.6}) is $C=1/c_\al$.

\vsk
\noindent
{\bf Conditioning of a least squares problem.} Suppose we are given $n$ complex numbers
$y_1, \ldots, y_n$ and we want to know whether there exists a polynomial $p$ of
degree at most $\al-1$ such that $p(j)=y_j$ for $1 \le j \le n$. Such a polynomial exists
if and only if
\begin{equation}
\de_k:=y_k-\left({\al \atop 1}\right)y_{k+1}+\left({\al \atop 2}\right)y_{k+2}
- \ldots + (-1)^\al y_{k+\al} = 0 \label{1.7}
\end{equation}
for $1 \le k \le n-\al$. Thus, to test the existence of $p$ we may compute
\[D(y_1, \ldots, y_n)=\left(\sum_{k=1}^{n-\al}\de_k^2\right)^{1/2}\]
and ask whether this
is a small number. Let $\cP_\al$ denote the set of all polynomials of degree at most $\al-1$
and put
\[E(y_1, \ldots, y_n)=\min_{p \in \cP_\al}\left(\sum_{j=1}^n |y_j-p(j)|^2\right)^{1/2}.\]
The question is whether $E(y_1, \ldots, y_n)$ may be large although $D(y_1, \ldots, y_n)$
is small. The answer to this question is (unfortunately)
in the affirmative and is in precise form given by the
formula
\begin{equation}
\max_{D(y_1, \ldots, y_n) \neq 0} \frac{E(y_1, \ldots, y_n)}{D(y_1, \ldots, y_n)}
\sim \frac{n^\al}{\sqrt{c_\al}}. \label{1.7a}
\end{equation}

\vsk
Here is our main result on the constants $c_\al$ we have encountered in the
five problems.

\vsk
\noindent
{\bf Theorem.} {\em We have the asymptotics
\begin{equation}
c_\al=\sqrt{8\pi \al}\,\left(\frac{4\al}{e}\right)^{2\al}\,
\left[1+O\left(\frac{1}{\sqrt{\al}}\right)\right]\quad \mbox{as}\quad \al \to \iy
\label{1.8}
\end{equation}
and the bounds}
\begin{equation}
\frac{4\al-2}{4\al^2-\al}\,\frac{(4\al)![\al !]^2}{[(2\al)!]^2}
\le c_\al \le \frac{4\al+1}{2\al+1}\,\frac{(4\al)![\al !]^2}{[(2\al)!]^2}
\quad \mbox{{\em for every}}\quad\al \ge 1.
\label{1.9}
\end{equation}

\vsk
In connection with (\ref{1.9}), notice that
\[\frac{(4\al)![\al !]^2}{[(2\al)!]^2} \sim \frac{1}{2}\,
\sqrt{8\pi \al}\,\left(\frac{4\al}{e}\right)^{2\al}.\]
Thus, the upper bound in (\ref{1.9}) is asymptotically exact, while the lower bound
in (\ref{1.9}) is asymptotically by the factor $1/(2\al)$ too small. This last
defect is nasty, but on the other hand it is clear that
$1/(2\al)$ is nothing in comparison with the astronomical growth of $(4\al/e)^{2\al}$.

\vsk
We discuss the five problems quoted here in more detail in Section \ref{S2}.
The theorem will be proved in Section \ref{S3}. Section \ref{S4} is devoted
to an alternative approach to Wirtinger-Sobolev inequalities and gives a new proof of the
coincidence of the constants in all the five problems.

\section{{\large Equivalence and history of the five problems}}\label{S2}

{\bf Toeplitz eigenvalues.} For $\al=1$, formula (\ref{1.1}) goes back to
Kac, Murdock, Szeg\"o \cite{KMS}.
In the late 1950's, Seymour Parter and the second of the authors
started tackling the general case, with Parter
embarking on the Toeplitz case and the second of us on the Wiener-Hopf case.
In \cite{Par99} ($\al=2$) and then in \cite{ParBAMS}, \cite{Par100} (general $\al$),
Parter established (\ref{1.1}).

\vsk
Subsequently, it turned out that the approach developed in
\cite{Wi88}, \cite{Wi100}, \cite{Wi106} can also be used to derive (\ref{1.1}).
This approach is as follows. Let $[T_n^\me(a)]_{j,k}$ be the $j,k$ entry
of $T_n^\me(a):=(T_n(a))^\me$ and consider the functions
\begin{equation}
n\,[T_n^\me(a)]_{[nx],[ny]}, \quad (x,y) \in [0,1]^2, \label{2.1}
\end{equation}
where $[nz]$ is the smallest integer in $\{1, \ldots, n\}$ that is greater than
or equal to $nz$. Let $K^\hn$ denote the integral operator on $L^2(0,1)$ with
the kernel (\ref{2.1}). One can prove two things. First,
\begin{equation}
\left\n \frac{1}{n^{2\al}}\,K^\hn - \frac{1}{b(1)}\,V_\al \right\n \to
0 \quad \mbox{as} \quad n \to \iy,
\label{2.2}
\end{equation}
where $V_\al$ is an integral operator with a certain completely identified kernel
$F_\al(x,y)$. And secondly, the eigenvalues of $K^\hn$ are just the eigenvalues of
$T_n^\me(a)$. These two insights imply that
\[\frac{1}{n^{2\al}}\,\frac{1}{\la_{\min}(T_n(a))} = \frac{1}{n^{2\al}}\,\la_{\max}(K^\hn)
\to \frac{1}{b(1)}\,\la_{\max}(V_\al)\]
or equivalently,
\[\la_{\min}(T_n(a)) \sim \frac{1/\la_{\max}(V_\al)}{n^{2\al}}\,b(1).\] 
The kernel $F_\al(x,y)$ is quite complicated, but it resembles
the kernel $G_\al(x,y)$ given by (\ref{1.5}).

\vsk
\noindent
{\bf Green's kernel.} In \cite{ParBAMS} and \cite{Wi100} it was further established
that $F_\al(x,y)$ is the Green kernel for the boundary problem (\ref{1.2}), (\ref{1.3}).
This implies at once that actually $F_\al(x,y)=G_\al(x,y)$ and $V_\al =K_\al$.
Thus, at this point is clear that in the first three problems of the introduction
we have to deal with one and the same constant $c_\al$.

\vsk
Expression (\ref{1.5}) was found in \cite{BoRS}. That paper concentrates on the case where
$b=1$, that is, where $a(t)=|1-t|^{2\al}$ ($t \in \bT$). Using a formula by Duduchava and Roch
for the inverse of $T_n(|1-t|^{2\al})$, it is shown in a direct way that
\[n^{1-2\al}\,[T_n^\me(|1-t|^{2\al})]_{[nx],[ny]} \to G_\al(x,y)
\quad \mbox{in} \quad L^\iy([0,1]^2).\]
Moreover, \cite{BoRS} has a short, self-contained, and elementary proof of the fact that
$G_\al(x,y)$ is the Green kernel of (\ref{1.2}), (\ref{1.3}).

\vsk
Rambour and Seghier \cite{RSIEOT} showed that
\begin{equation}
n^{1-2\al}\,[T_n^\me(|1-t|^{2\al}b(t))]_{[nx],[ny]} \to \frac{1}{b(1)}\,G_\al(x,y)
\quad \mbox{in} \quad L^\iy([0,1]^2) \label{2.3a}
\end{equation}
under the assumptions on $b$ made in the introduction. Evidently, (\ref{2.3a})
implies (\ref{2.2}) (but not vice versa). For $\al=1$,
result (\ref{2.3a}) was known from previous work of
Courant, Friedrichs, and Lewy \cite{CFL} and Spitzer and Stone \cite{SpiSto}.
The authors of \cite{RSIEOT} were obviously not aware of papers \cite{ParBAMS}
and \cite{Wi100} and rediscovered again that $G_\al(x,y)$ is Green's kernel of (\ref{1.2}),
(\ref{1.3}).

\vsk
\noindent
{\bf Wirtinger-Sobolev.} The connection between the minimal eigenvalue of
(\ref{1.2}), (\ref{1.3}) and the best constant in (\ref{1.6}) with the boundary
conditionÇs (\ref{1.3}) is nearly obvious. Indeed, we have
\[c_\al = \min\,\frac{((-1)^\al u^{(2\al)},u)}{(u,u)},\]
where $(\cdot,\cdot)$ is the inner product in $L^2(0,1)$ and the minimum is over
all nonzero and smooth functions $u$ satisfying (\ref{1.3}). Upon $\al$ times
partially integrating and using the boundary conditions, one gets
\[c_\al=\min\,\frac{(u^\hal, u^\hal)}{(u,u)}
= \min\,\frac{\int_0^1|u^\hal(x)|^2dx}{\int_0^1|u(x)|^2dx},\]
which is equivalent to saying that the best constant $C$ in (\ref{1.6})
with the boundary conditions (\ref{1.3}) is $C=1/c_\al$.

\vsk
Numerous versions of inequalities of the Wirtinger-Sobolev type have been
established for many decades. The original inequality says that
\begin{equation}
\int_0^1|u(x)|^2dx -\left|\int_0^1 u(x)dx\right|^2
\le \frac{1}{(2\pi)^2}\int_0^1|u'(x)|^2dx \label{2.4}
\end{equation}
whenever $u \in C^1[0,1]$ and $u(0)=u(1)$. This inequality appears in
different modifications, sometimes with the additional requirement that
$\int_0^1 u(x)dx=0$ and frequently over the interval $(0,2\pi)$, in which
case the constant $1/(2\pi)^2$ becomes $1$
(see, e.g., \cite[pp. 184-187]{HLP}). The proof of (\ref{2.4}) is in fact very
simple: take the Fourier expansion $u(x)=\sum u_k e^{2\pi i k x}$ and use
Parseval's equality. We will say more on this topic in Section~\ref{S4},
which contains a direct proof of the fact that the best constant $C$ in
(\ref{1.6}) with the boundary conditions (\ref{1.3}) is the inverse of
the constant $c_\al$ of (\ref{1.1}).

\vsk
\noindent
{\bf The least squares problem.} The least squares result is from \cite{BoTest}.
Define the linear operator $\nabla : \bC^n \to \bC^n$ by
$\nabla (y_1, \ldots, y_n)=(\de_1, \ldots, \de_{n-\al}, 0, \ldots, 0)$,
where the $\de_k$'s are given by (\ref{1.7}), put ${\rm Ker}\,\nabla
:=\{y \in \bC^n: \nabla y =0\}$, and denote by $P_{{\rm Ker}\,\nabla}$
the orthogonal projection of $\bC^n$ onto ${\rm Ker}\,\nabla$.
The left-hand side of (\ref{1.7a})
is nothing but
\begin{equation}
\max_{y \notin {\rm Ker}\,\nabla}\,\frac{\n y- P_{{\rm Ker}\nabla}\, y\n_2}{\n\nabla y\n_2},
\label{2.9}
\end{equation}
where $\n \cdot \n_2$ is the $\ell^2$ norm on $\bC^n$. With $\nabla^+$ denoting the
Moore-Penrose inverse of $\nabla$, we have the equality $I - P_{{\rm Ker}\,\nabla}
=\nabla^+\nabla$. This shows that (\ref{2.9}) is the norm of $\nabla^+$, that is,
the inverse of the smallest nonzero singular value of $\nabla$. But $\nabla \nabla^*$
can be shown to be of the form
\[J\left(\begin{array}{cc}
T_{n-\al}(|1-t|^{2\al}) & 0\\
0 & O_\al \end{array}\right)J,\]
where $J$ is a permutation matrix and $O_\al$ is the $\al \times \al$ zero matrix.
Thus, the smallest nonzero singular value of $\nabla$ is the square root of
$\la_{\min}(T_{n-\al}(|1-t|^{2\al})) \sim c_\al/n^{2\al}$, which brings us back to the
beginning.

\vsk
\noindent
{\bf A wrong conjecture.} The first three values of $c_\al$ are
\[c_1=\pi^2=9.8696, \quad c_2=500.5467, \quad c_3= 61529,\]
and the first three values of $((\al+1)\pi/2)^{2\al}$ are
\[\pi^2=9.8696, \quad 493.1335, \quad 61529.\]
We all know that one should not guess the asymptotics of a sequence
from its first three terms. But because of the amazing coincidence in the
case $\al=3$, it is indeed tempting to conjecture that $c_\al \sim ((\al+1)\pi/2)^{2\al}$.
Our main result shows that this conjecture is wrong. The first three values of the
correct asymptotics $c_\al \sim \sqrt{8 \pi \al}\,(4\al/e)^{2\al}$ are
\[10.8555, \quad 531.8840, \quad 64269.\]

\section{{\large Proof of the main result}}\label{S3}

We employ the equality $1/c_\al=\n K_\al\n$, where $K_\al$ is the integral operator
on $L^2(0,1)$ with kernel (\ref{1.5}).

\vsk
\noindent
{\bf The kernel's peak.} It will turn out that the main contribution to the kernel
$G_\al(x,y)$ comes from a neighborhood of $(\frac{1}{2}, \frac{1}{2})$, and so for
later convenience we consider instead the integral operator $\widetilde{K}_\al$
on $L^2(-1,1)$ whose kernel is
\[\widetilde{G}_\al(x,y)=\frac{1}{2}\,G_\al \left(\frac{1+x}{2}, \frac{1+y}{2}\right).\]
The operator $\widetilde{K}_\al$ has the same norm as $K_\al$,
its kernel is symmetric about $(0,0)$, and the main contribution to the kernel
comes from a neighborhood of $(0,0)$. If we make the substitution $t \to (1+t)/2$
in the integral we see that
\[\widetilde{G}_\al(x,y)=\frac{1}{[(\al-1)!]^2}\,\frac{1}{4^\al}\,H_\al(x,y)\]
with
\begin{equation}
H_\al(x,y)=(1+x)^\al (1+y)^\al \int_{\max(x,y)}^1
\frac{(t-x)^{\al-1}(t-y)^{\al-1}}{((1+t)/2)^{2\al}}\,dt \label{3.1}
\end{equation}
when $x+y \ge 0$. We shall show that $H_\al(x,y)$ is equal to
$(1/\al)(1-x^2)^\al (1-y^2)^\al$ plus a kernel whose norm is
smaller by a factor $O(1/\sqrt{\al})$.

\vsk
The logarithmic derivative in $t$ of the function
$(t-x)\,(t-y)/((1+t)/2)^2$
is
\[\frac{(2+x+y)t-x-y-2xy}{(1+t)(t-x)(t-y)},\]
which is positive for $t > \max (x,y)$. (Recall that we are in the case
$x+y \ge 0$.) Hence the function achieves its maximum $(1-x)(1-y)$
at $t=1$ and nowhere else. The function
$(1+x)(1+y)(1-x)(1-y)$ achieves its maximum at $x=0$, $y=0$ and nowhere else.
Putting these together we see that the function
\[(1+x)(1+y)\,\frac{(t-x)(t-y)}{((1+t)/2)^2}\]
achieves its maximum $1$ at $t=1$, $x=0$, $y=0$, and outside a neighborhood
of this point, say outside the set
$t \ge 1-\eps$, $|x| \le \eps$, $|y| \le \eps$, there is a bound
\[(1+x)(1+y)\,\frac{(t-x)(t-y)}{((1+t)/2)^2} < 1-\de\]
for some $\de >0$. It follows that outside the same neighborhood the
integrand in (\ref{3.1}) with its outside factor is $O((1-\de)^\al)$.
This is also the bound after we integrate. We take any $\eps < 1/2$,
and have shown that
\[H_\al(x,y) =(1+x)^\al (1+y)^\al \chi_\eps(x)\chi_\eps(y)
\int_{1-\eps}^1 \frac{(t-x)^{\al-1}(t-y)^{\al-1}}{((1+t)/2)^{2\al}}
+O((1-\de)^\al),\]
where $\chi_\eps$ is $1$ on $[-\eps,\eps]$ and zero elsewhere. Substituting
$t=1-\tau$ we arrive at the formula
\begin{eqnarray}
H_\al(x,y)& = & (1+x)^\al (1+y)^\al \chi_\eps(x)\chi_\eps(y)
\int_0^\eps \left[\frac{\left(1-\frac{\tau}{1-x}\right)
\left(1-\frac{\tau}{1-y}\right)}{\left(1-\frac{\tau}{2}\right)^2}
\right]^\al \times \nonumber\\
& & \hspace{20mm} \times \,\frac{d\tau}{\left(1-\frac{\tau}{1-x}\right)
\left(1-\frac{\tau}{1-y}\right)} +O((1-\de)^\al). \label{3.2}
\end{eqnarray}

\vsk
\noindent
{\bf The kernel's asymptotics.} Let us compute the asymptotics of the kernel.
The choice $\eps < 1/2$ guarantees that $\tau/(1-x)$, $\tau/(1-y)$, $\tau/2$
belong to $(0,1)$. This implies that
\[\frac{\left(1-\frac{\tau}{1-x}\right)
\left(1-\frac{\tau}{1-y}\right)}{\left(1-\frac{\tau}{2}\right)^2}
=e^{-\tau \ph(x,y)+O(\tau^2)}\]
with \[\ph(x,y)=\frac{1-xy}{(1-x)(1-y)}.\]
We split the integral in (\ref{3.2}) into $\int_0^{1/\sqrt{\al}}$ and
$\int_{1/\sqrt{\al}}^\eps$. If $\eps >0$ is small enough, which we may assume,
the term $O(\tau^2)$ is at most $\tau \ph(x,y)/2$ in absolute value.
Hence the integral $\int_{1/\sqrt{\al}}^\eps$ is at most
\[\int_{1/\sqrt{\al}}^\eps e^{-\al \tau \ph(x,y)/2} O(1) d\tau=O\left(e^{-\ga_1
\sqrt{\al}}\right)\]
with some $\ga_1>0$. For $\tau < 1/\sqrt{\al}$ we have $\al \tau^2 <1$ and hence
$e^{\al O(\tau^2)}=1+\al O(\tau^2)$. Consequently, the integral
$\int_0^{1/\sqrt{\al}}$  is equal to
\begin{eqnarray}
& & \int_0^{1/\sqrt{\al}} e^{-\al \tau \ph(x,y)}e^{\al O(\tau^2)} (1+O(\tau))d\tau
\nonumber\\
& & = \int_0^{1/\sqrt{\al}} e^{-\al \tau \ph(x,y)} \left(1+O(\tau)+\al O(\tau^2)
\right)d\tau.
\label{3.3}
\end{eqnarray}
Since, for $k=0,1,2$,
\[\int_{1/\sqrt{\al}}^\iy \tau^ke^{-\al \tau \ph(x,y)}d\tau
=O\left(e^{-\ga_2 \sqrt{\al}}\right)\]
with $\ga_2 >0$ and
\[\int_0^\iy \tau^ke^{-\al \tau \ph(x,y)}d\tau=O\left(\frac{1}{\al^{k+1}}\right),\]
it follows that (\ref{3.3}) is
\begin{eqnarray*}
& & \int_0^\iy e^{-\al \tau \ph(x,y)}\left(1+O(\tau)+\al O(\tau^2)\right)d\tau
+O\left(e^{-\ga_2 \sqrt{\al}}\right)\\
& & =\left.\frac{e^{-\al \tau \ph(x,y)}}{-\al \ph(x,y)}\right|_0^\iy
+O\left(\frac{1}{\al^2}\right)+\al\,O\left(\frac{1}{\al^3}\right)
+O\left(e^{-\ga_2\sqrt{\al}}\right)\\
& & = \frac{1}{\al \ph(x,y)} +O\left(\frac{1}{\al^2}\right).
\end{eqnarray*}
In summary,
\[H_\al(x,y)=\frac{(1-x^2)^\al(1-y^2)^\al}{\al \ph(x,y)}
\chi_\eps(x)\chi_\eps(y)+O\left(\frac{1}{\al^2}\right),\]
uniformly for $|x|, |y| \le \eps$. Expanding near $x=y=0$ we obtain
\begin{equation}
H_\al(x,y)=\frac{1}{\al}\,(1-x^2)^\al(1-y^2)^\al (1+O(x)+O(y))+
O\left(\frac{1}{\al^2}\right),
\label{3.4}
\end{equation}
again uniformly. This was derived for $|x|, |y| \le \eps$, but because of
(\ref{3.2}) we see that this holds uniformly for all $x$ and $y$ satisfying $x+y \ge 0$.
This last condition can also be dropped by the symmetry of $H_\al(x,y)$.

\vsk
\noindent
{\bf The asymptotics of the norm.} If an integral operator $K$ is of the form
\[(Ku)(x)=\int_{-1}^1 f(x)g(y)u(y)dy,\]
then $\n K\n =\n f\n_2 \n g\n_2$, where $\n \cdot\n_2$ is the norm in $L^2(-1,1)$.
Let us denote the integral operator with the kernel $H_\al(x,y)$ by $M_\al$.
Furthermore, in view of (\ref{3.4}) we denote by $M_\al^0, M_\al^1, M_\al^2$
the integral operators with the kernels
\[(1-x^2)^\al(1-y^2)^\al, \quad O(x)\,(1-x^2)^\al(1-y^2)^\al, \quad
O(y)\,(1-x^2)^\al(1-y^2)^\al,\]
respectively. {From} (\ref{3.4}) we infer that
\[\n M_\al\n =\frac{1}{\al}\,\n M_\al^0 +M_\al^1
+M_\al^2\n+O\left(\frac{1}{\al^2}\right).\]
Since
\begin{eqnarray*}
\n M_\al^0\n & = & \int_{-1}^1(1-x^2)^{2\al}dx =\sqrt{\frac{\pi}{2\al}}
\left(1+O\left(\frac{1}{\al}\right)\right),\\
\n M_\al^1\n^2 & = & \int_{-1}^1 O(x^2)\,(1-x^2)^{2\al}dx
\int_{-1}^1(1-y^2)^{2\al}dy\\
& = & O\left(\frac{1}{4\al}\,\sqrt{\frac{\pi}{2\al}}\right)\,
\sqrt{\frac{\pi}{2\al}}= O\left(\frac{1}{\al^2}\right),
\end{eqnarray*}
and a similar estimate is valid for $\n M_\al^2\n^2$, we finally get
\begin{eqnarray*}
\n K_\al\n & = & \n \widetilde{K}_\al\n = \frac{1}{[(\al-1)!]^2}\,\frac{1}{4^{2\al}}
\,\n M_\al\n\\
& = & \frac{1}{[(\al-1)!]^2}\,\frac{1}{4^{2\al}\,\al}
\left[\n M_\al^0\n +O\left(\n M_\al^1\n\right) +O\left(\n M_\al^2\n\right)
+O\left(\frac{1}{\al}\right)\right]\\
& = & \frac{\al^2}{\al^{2\al}e^{-2\al}2\pi \al}\,\left(1+O\left(\frac{1}{\al}\right)\right)
\frac{1}{4^\al\,\al}\,\left[\sqrt{\frac{\pi}{2\al}}+O\left(\frac{1}{\al}\right)\right]\\
& = & \left(\frac{e}{4\al}\right)^{2\al}\,\frac{1}{\sqrt{8\pi \al}}\,
\left(1+O\left(\frac{1}{\sqrt{\al}}\right)\right),
\end{eqnarray*}
which is the same as (\ref{1.8}).

\vsk
\noindent
{\bf The lower bound.} To prove the lower bound in (\ref{1.9}), we start with (\ref{3.1})
and the inequality
\[\frac{(t-x)(t-y)}{((1+t)/2)^2}\le (1-x)(1-y),\]
which was established in the course of the above proof.
If $x+y \ge 0$, then $\max(x,y) \ge 0$ and consequently,
\begin{eqnarray*}
H_\al(x,y) & \le &  (1+x)^\al (1+y)^\al \int_0^1 (1-x)^{\al-1}(1-y)^{\al-1}\,
\frac{dt}{((1+t)/2)^2}\\
& = & 2(1+x)(1-x^2)^{\al-1}(1+y)(1-y^2)^{\al-1}.
\end{eqnarray*}
Hence
\[\n M_\al\n \le 2\int_{-1}^1(1+x)^2(1-x^2)^{2\al-2}dx
=4^{2\al}\,\frac{(2\al)!(2\al-2)!}{(4\al-1)!}\]
and thus
\begin{eqnarray*}
\n K_\al\n & = & \n \widetilde{K}_\al \n
\le \frac{1}{[(\al-1)!]^2}\,\frac{1}{4^{2\al}}\,4^{2\al}\,\frac{(2\al)!(2\al-2)!}{(4\al-1)!}\\
& = & \frac{\al^2}{\al!\al!}\,\frac{(2\al)!(2\al)!(4\al-1)}{(2\al-1)(2\al)(4\al)!},
\end{eqnarray*}
which is equivalent to the assertion.

\vsk
\noindent
{\bf The upper bound.} The proof of the upper bound in (\ref{1.9})
is based on the observation that $1/c_\al$
is the best constant for which the inequality
\[\int_0^1|u(x)|^2dx \le \frac{1}{c_\al}\,\int_0^1|u^\hal (x)|^2 dx\]
is true for all $u \in C^{\al}[0,1]$ satisfying $u^{(j)}(0)=u^{(j)}(1)=0$
for $0 \le j \le \al-1$. If we insert $u(x)=x^\al(1-x)^\al$, the inequality
becomes
\[\int_0^1x^{2\al}(1-x)^{2\al}dx \le \frac{1}{c_\al}\,\int_0^1
\left[\frac{d^\al}{dx^\al}\left(x^\al(1-x)^\al\right)\right]^2\,dx.\]
The integral on the left is $[(2\al)!]^2/(4\al+1)!$, and in the integral
on the right we make the substitution $x=(1+y)/2$ to get
\[\int_0^1
\left[\frac{d^\al}{dx^\al}\left(x^\al(1-x)^\al\right)\right]^2\,dx
=
\frac{1}{4^\al}\,\int_{-1}^1 \left[\frac{d^\al}{dy^\al}\left(y^2-1\right)^\al
\right]^2\frac{dy}{2}.\]
The function $(d^\al/dy^\al)(y^2-1)^\al$ is $2^\al\,\al!$ times the usual Legendre
polynomial $P_\al(y)$ and it is well known
that $\n P_\al\n_2^2=2/(2\al+1)$ (see, for example, \cite{JEL}).
Consequently, the integral on the right is
\[\frac{1}{4^\al}\,\frac{1}{2}\,2^{2\al}(\al!)^2\,\frac{2}{2\al+1}
=\frac{(\al!)^2}{2\al+1}.\]
In summary,
\[\frac{[(2\al)!]^2}{(4\al+1)!} \le \frac{1}{c_\al}\,\frac{(\al!)^2}{2\al+1},\]
which is the asserted inequality.

\vsk
\noindent
{\bf Refinements.} By carrying out the approximations further we could refine
(\ref{3.2}) to the form
\[H_\al(x,y)=\frac{1}{\al}\,(1-x^2)^\al(1-y^2)^\al\left(1+\sum_{i \ge 1, j \ge 0}
\frac{p_{ij}(x,y)}{\al^i}\right),\]
where each $p_{ij}(x,y)$ is a homogeneous polynomial of degree $j$. The operator with kernel
$(1-x^2)^\al(1-y^2)^\al p_{ij}(x,y)$ has norm of the order $\al^{-(j+1)/2}$, so we get
further approximations to $H_\al$ in this way, whence further approximations to the norm.
(We do this by using the fact that the nonzero eigenvalues of a finite-rank kernel $\sum_{i=1}^mf_i(x)\,g_i(y)$ are the same as those of the $m\times m$ matrix whose $i,j$ entry is the inner product $(f_i,\,g_j)$. One can see from this in particular that, because of evenness and oddness, with each approximation the power of $\al$ goes down by one.) However, these would probably not be of great interest.

\section{{\large Another approach to Wirtinger-Sobolev inequalities}}\label{S4}

We now show how Wirtinger-Sobolev integral inequalities can be derived from their
discrete analogues, which, in dependence on the boundary conditions, are inequalities
for circulant or Toeplitz matrices. In the Toeplitz case, we get in this way a new proof
of the fact that the constants in the first and fourth problems are the same.

\vsk
Discrete versions of Wirtinger-Sobolev type
inequalities were first established by Fan, Taussky, and Todd \cite{FTT}, and the
subject has been developed further since then (see, for example, \cite{MM}
and the references therein). In particular, for circulant matrices
the following is not terribly new, but
it fits very well with the topic of this paper and perfectly illustrates
the difference between the circulant and Toeplitz cases.

\vsk
\noindent
{\bf Circulant matrices.} For a Laurent polynomial $a(t)=\sum_{k=-r}^r a_kt^k$
($t \in \bT$) and $n \ge 2r+1$, let $C_n(a)$ be the $n \times n$ circulant matrix whose first
row is
\[(a_0, a_{-1}, \ldots, a_{-r}, 0, \ldots, 0, a_r, a_{r-1}, \ldots, a_1).\]
Thus, $C_n(a)$ results from the Toeplitz matrix $T_n(a)$ by periodization.
The singular values of $C_n(a)$ are $|a(\om_n^j)|$ ($j=1, \ldots, n$), where
$\om_n=e^{2\pi i/n}$.

\vsk
Now let $a(t)=(1-t)^\al$ ($t \in \bT$). One of the singular values of
$C_n(a)$ is zero, which causes a slight complication. It is easily seen that
${\rm Ker}\,C_n(a)={\rm span}\,\{(1,1, \ldots,1)\}$. With notation as in Section \ref{S2},
$I-P_{{\rm Ker}\,C_n(a)}=C_n^+(a)C_n(a)$ and hence
\begin{equation}
\n u-P_{{\rm Ker}\,C_n(a)}u\n_2 \le \n C_n^+(a)\n\,\n C_n(a)u\n_2 \label{4.1}
\end{equation}
for all $u$ in $\bC^n$ with the $\ell^2$ norm. The inverse of the (spectral) norm of the
Moore-Penrose inverse $C_n^+(a)$ is the smallest nonzero singular
value of $C_n(a)$ and consequently,
\begin{equation}
\frac{1}{\n C_n^+(a)\n}=
|1-\om_n|^\al =\left(4\,\sin^2\frac{\pi}{n}\right)^{\al/2} \sim \frac{(2\pi)^\al}{n^\al}
\label{4.2}
\end{equation}
The projection $P_{{\rm Ker}\,C_n(a)}$ acts by the rule
\begin{equation}
P_{{\rm Ker}\,C_n(a)}u=\left(\frac{1}{n}\sum_{j=1}^n u_j, \ldots,
\frac{1}{n}\sum_{j=1}^n u_j\right). \label{4.3}
\end{equation}
Inserting (\ref{4.2}) and (\ref{4.3}) in (\ref{4.1}) we get
\begin{eqnarray}
\n C_n(a)u\n^2_2 & \ge & \left(4\,\sin^2\frac{\pi}{n}\right)^{\al}
\sum_{i=1}^n\left|u_i-\frac{1}{n}\sum_{j=1}^n u_j\right|^2\nonumber\\
& = & \left(4\,\sin^2\frac{\pi}{n}\right)^{\al} \left(\sum_{i=1}^n |u_i|^2
-\frac{1}{n}\left|\sum_{i=1}^n u_i \right|^2\,\right). \label{4.4}
\end{eqnarray}
This is called a (higher-order) discrete Wirtinger-Sobolev inequality
and was by different methods established in \cite{MM}.

\vsk
\noindent
{\bf Periodic boundary conditions.} As already said, the wanted inequality (\ref{4.5})
follows almost immediately from Parseval's identity. So the following might seem
unduly complicated. However, the analogue of (\ref{4.5}) for zero boundary conditions
is not straightforward from Parseval's identity, whereas just the following
also works in that case.

\vsk
Let $u$ be a $1$-periodic function in
$C^\iy(\bR)$. We apply (\ref{4.4}) to $u_n=(u(j/n))_{j=1}^n$. The $j$th component of
$C_n(a)u_n$ is
\[u\left(\frac{j}{n}\right)-\left({\al \atop 1}\right)u\left(\frac{j+1}{n}\right)
+ \ldots + (-1)^\al u\left(\frac{j+\al}{n}\right)
= u^\hal\left(\frac{j}{n}\right)\,\frac{1}{n^\al}+O\left(\frac{1}{n^{\al+1}}\right),\]
the $O$ being independent of $j$. It follows that
\[\n C_n(a)u\n^2_2 =\left(\sum_{j=1}^n \left|u^\hal \left(\frac{j}{n}\right)\right|^2\,
\frac{1}{n^{2\al}}\,\right)+O\left(\frac{1}{n^{2\al}}\right).\]
Consequently, multiplying (\ref{4.4}) by $n^{2\al-1}$ and passing to the limit $n \to \iy$ we
arrive at the inequality
\begin{equation}
\int_0^1|u^\hal(x)|^2dx \ge (2\pi)^{2\al}\left(\int_0^1|u(x)|^2dx
-\left|\int_0^1 u(x)dx\right|^2\,\right). \label{4.5}
\end{equation}

Assume finally that $u\in C^\al[0,1]$ and $u^{(j)}(0)=u^{(j)}(1)$ for $0 \le j \le \al-1$.
We have $u(x)=\sum_{k=-\iy}^\iy e^{2\pi i k x}$ with
\[u_k=\int_0^1u(x)e^{-2 \pi i k x} dx.\]
We integrate the last equality $\al$ times partially and use the boundary conditions
to obtain that
\[|u_k| =\frac{1}{(2\pi |k|)^{\al}}\,\left|\int_0^1 u^\hal(x)e^{-2\pi ikx}dx\right|.\]
Since $u^\hal \in L^2(0,1)$, we see that $|u_k|=v_k \,O(1/|k|^\al)$ with $\sum_{k=-\iy}^\iy
v_k^2 < \iy$. This implies that
\begin{equation}
\sum_{k=-\iy}^\iy |k|^{2\al} |u_k|^2 < \iy. \label{4.6}
\end{equation}

We know that (\ref{4.5}) is true with $u(x)$ replaced by $(S_Nu)(x)=\sum_{k=-N}^N
u_ke^{2\pi i k x}$,
\begin{equation}
\int_0^1 |(S_Nu)^\hal(x)|^2dx \ge (2\pi)^{2\al}\left(\int_0^1|(S_Nu)(x)|^2dx
-\left|\int_0^1 (S_Nu)(x)dx\right|^2\,\right). \label{4.7}
\end{equation}
{From} (\ref{4.6}) we infer that
\begin{eqnarray*}
& & \int_0^1|u^\hal(x)|^2dx-\int_0^1 |(S_Nu)^\hal(x)|^2dx
=\sum_{|k|>N}|k|^{2\al}|u_k|^2=o(1),\\
& & \int_0^1|u(x)|^2dx-\int_0^1 |(S_Nu)(x)|^2dx
=\sum_{|k|>N}|u_k|^2=o(1),
\end{eqnarray*}
and since $\int_0^1(S_Nu)(x)dx=\int_0^1 u(x)dx =u_0$, passage to the limit $N \to \iy$
in (\ref{4.7}) yields (\ref{4.5}) under the above assumptions on $u$.

\vsk
\noindent
{\bf Toeplitz matrices.} Again let $a(t)=(1-t)^\al$ ($t \in \bT$), but consider now
the Toeplitz matrix $T_n(a)$ instead the circulant matrix $C_n(a)$. It can be easily
verified or deduced from \cite[formula~(2.13)]{BSUni}
or \cite[formula~(1.4)]{WiII} that
\[T_n^*(a)T_n(a)=T_n(b)-R_{\al}\]
where $b(t)=|1-t|^{2\al}$ and $R_\al$ is a matrix of the form
\[R_\al=\left(\begin{array}{cc}S_\al & 0\\ 0 & O_{n-\al}\end{array}\right)\]
with an $\al \times \al$ matrix $S_\al$ independent of $n$.
Consequently,
\[\n T_n(a)u\n_2^2=(T_n(a)u,T_n(a)u)=(T_n(b)u,u)-(R_\al u,u).\]
It follows that
\begin{equation}
\n T_n(a)u\n_2^2 \ge \la_{\min}(T_n(b))\,\n u\n_2^2 -(R_\al u,u) \label{4.8}
\end{equation}
for all $u \in \bC^n$. This is the Toeplitz analogue of (\ref{4.4}).

\vsk
\noindent
{\bf Zero boundary conditions.} Let $u \in C^\iy(\bR)$ be a function which
vanishes identically outside $(0,1)$. As in the circulant case, we replace the
$u$ in (\ref{4.8}) by $u_n=(u(j/n))_{j=1}^n$, multiply the result by $n^{2\al-1}$
and pass to the limit $n \to \iy$. Taking into account that $\la_{\min}(T_n(b))
\sim c_\al/n^{2\al}$, we obtain
\begin{equation}
\int_0^1|u^\hal(x)|^2dx \ge c_\al\int_0^1|u(x)|^2-\lim_{n\to \iy}
n^{2\al-1}(R_\al u_n,u_n). \label{4.9}
\end{equation}
By assumption, $u$ and all its derivatives vanish at $0$. This implies that
\[u\left(\frac{j}{n}\right)=\sum_{k=0}^{\al-1}\frac{u^{(k)}(0)}{k!}\,\frac{j^k}{n^k}
+\frac{u^{(\al)}(\xi_{j,n})}{\al!}\,\frac{j^{\al}}{n^{\al}}
=O\left(\frac{1}{n^{\al}}\right)\]
for each fixed $j$. Since $(R_\al u_n,u_n)$ is a bilinear form of $u(1/n), \ldots,
u(\al/n)$, we arrive at the conclusion that $(R_\al u_n,u_n)=O(1/n^{2\al})$.
Hence, (\ref{4.9}) is actually the desired inequality
\begin{equation}
\int_0^1 |u^\hal (x)|^2 dx \ge c_\al \int_0^1 |u(x)|^2 dx. \label{4.10}
\end{equation}

\vsk
The approximation argument employed in the case of periodic boundary conditions
is also applicable in the case at hand and allows us to relax the $C^\iy$
assumption. It results that (\ref{4.10}) is valid for every $u \in C^\al[0,1]$
satisfying $u^{(j)}(0)=u^{(j)}(1)=0$ for $0 \le j \le \al-1$.

\vsg
\noindent
\begin{minipage}[t]{8cm}
Albrecht B\"ottcher\\
Fakult\"at f\"ur Mathematik\\
TU Chemnitz\\
09107 Chemnitz\\
Germany\\[1ex]
aboettch@mathematik.tu-chemnitz.de
\end{minipage}
\begin{minipage}[t]{6.5cm}
Harold Widom \\
Department of Mathematics\\
University of California\\
Santa Cruz, CA 95064\\
USA\\[1ex]
widom@math.ucsc.edu

\end{minipage}

\end{document}